\documentclass{amsart}
\usepackage{amsmath,amsthm,amsfonts,amssymb}
\usepackage{mathrsfs}
\usepackage{enumitem}
\usepackage[colorlinks,linkcolor=blue,anchorcolor=blue,citecolor=blue,urlcolor=blue]{hyperref}
\usepackage[initials,nobysame]{amsrefs}

\renewcommand{\PrintDOI}[1]{DOI \href{https://doi.org/#1}{\nolinkurl{#1}}}

\makeatletter
\@namedef{subjclassname@2020}{\textup{2020} Mathematics Subject Classification}
\makeatother
\newtheorem{thm}{Theorem}[section]
\newtheorem{lem}[thm]{Lemma}
\newtheorem{prop}[thm]{Proposition}
\newtheorem{cor}[thm]{Corollary}
\theoremstyle{definition}
\newtheorem{dfn}[thm]{Definition}
\theoremstyle{remark}
\newtheorem{rem}[thm]{Remark}
\DeclareMathOperator{\Var}{Var}
\DeclareMathOperator{\osc}{osc}
\DeclareMathOperator{\diam}{diam}
\newcommand{\R}{\mathbb R}
\newcommand{\1}{\mathbf 1}
\newcommand{\C}{\mathcal C}

\allowdisplaybreaks[1]
\numberwithin{equation}{section}

\hypersetup{pdftitle={Variational Estimates for the Fractional Maximal Operator on Finite Graphs},pdfauthor={Shengzhe Lu and Fayou Zhao}}
\begin{document}
\title[Fractional variation on finite graphs]{Variational Estimates for the Fractional Maximal Operator on Finite Graphs}
\author[S. Lu and F. Zhao]{Shengzhe Lu and Fayou Zhao$^*$}
\address[Shengzhe Lu]{Department of Mathematics, Shanghai University, Shanghai 200444, P.R. China}
\email{lushengzhe@shu.edu.cn}
\address[Fayou Zhao]{Newtouch Center for Mathematics of Shanghai University, Department of Mathematics, Shanghai University, Shanghai 200444, P.R. China}
\email{fyzhao@shu.edu.cn}
\thanks{$\ast$ Corresponding author}
\date{}
\begin{abstract}
We study edge variation estimates for the centered fractional maximal operator on finite connected graphs.
For $0<\alpha<1$, we determine the optimal coefficient in the uniform estimate of the oscillation of the maximal function by the edge $p$-variation of the original function.
We also obtain two-sided variation bounds involving the number of vertices and the diameter.
On complete graphs, we prove that the sharp $(p,q)$ constant is
$(1-n^{\alpha-1})(n-1)^{1/q-1/p}$ for all $0<p\leq q\leq\infty$.
On star graphs with $n\geq4$, we determine the sharp infinite-exponent constant and show that the sharp constant for $p=2$ equals $1-n^{\alpha-1}$ if and only if $n^{\alpha-1}\geq1/2$.
Combining this with the known estimate for $p=1$, we also obtain sharp mixed estimates.
Extremizers and counterexamples are given.
\end{abstract}
\maketitle

\section{Introduction}\label{sec:intro}
The regularity of maximal functions has received considerable attention since Kinnunen \cite{Kinnunen} proved that the Hardy--Littlewood maximal operator is bounded on $W^{1,p}(\mathbb R^d)$ for $1<p\leq\infty$.
On a finite graph, differences along edges give a discrete counterpart of derivatives.
In this paper, we study the best constant in
\[
\Var_q(M_{\alpha,G}f)\leq\C_{\alpha,G;p,q}\Var_p(f),
\qquad 0<\alpha<1,\quad 0<p\leq q\leq\infty,
\]
where $M_{\alpha,G}$ is the centered fractional maximal operator with respect to vertex counting measure, and each undirected edge is counted once in the variation.

Liu and Xue \cite{LX} studied variational estimates for maximal functions on finite graphs.
Gonz\'alez-Riquelme and Madrid \cite{GRM1} obtained further results on complete graphs, star graphs, and general finite graphs; their results include mixed-variation boundedness for the fractional maximal operator.
They subsequently determined the sharp $p=2$ constant for the ordinary maximal operator on star graphs \cite{GRM2}.
Gonz\'alez-Riquelme, Kova\v{c}, and Madrid \cite{GKM} determined the sharp constant for the ordinary maximal operator on complete graphs for every finite positive exponent.
Liu, Lian, and Wang \cite{LLW} obtained sharp fractional variation estimates on complete and star graphs for $p=1$ and $q=1/(1-\alpha)$.

Our first result is the optimal uniform estimate
\[
\osc(M_{\alpha,G}f)\leq\Omega_{\alpha,p}(n)\Var_p(f)
\]
over all connected graphs with $n$ vertices; see Theorem~\ref{thm:sharp-osc-variation}.
Ordered vertex values give the upper bound, and functions on a path attain equality.
Combining this estimate with point-mass tests and an edge-count bound in terms of the diameter, we obtain two-sided bounds for an individual graph in Theorems~\ref{thm:general-pp} and~\ref{thm:generalbounds}.

For complete graphs, the maximal operator truncates all vertex values at the same level.
Starting from the mean-truncation inequality of \cite{GKM}*{Theorem 1}, we prove monotonicity of a normalized sum of truncated differences.
This gives the exact formula
\[
\C_{\alpha,K_n;p,q}=(1-n^{\alpha-1})(n-1)^{1/q-1/p},
\qquad n\geq2,\quad 0<p\leq q\leq\infty,
\]
with equality for every point mass; see Theorem~\ref{thm:complete-mixed}.

On a star graph, two-point balls contribute in addition to singletons and the whole graph.
For $n\geq4$, the sharp constant for $p=1$ is $1-n^{\alpha-1}$, as proved in \cite{LLW}*{proof of Theorem 1.6, (4.4)}.
The $p=2$ sharp constant equals this value precisely when $n^{\alpha-1}\geq1/2$; see Theorems~\ref{thm:fracstar-p1-sharp} and~\ref{thm:fracstar-p2-sharp}.
Theorem~\ref{thm:fracstarinf} gives the sharp infinite-exponent constant, and Theorem~\ref{thm:fracstar-mixed-sharp} gives the sharp mixed constants when $p\leq1\leq q$, or when $n^{\alpha-1}\geq1/2$ and $p\leq2\leq q$.
Proposition~\ref{prop:fracstar-counter} shows that the equal-exponent formulas do not extend to all smaller or larger exponents.

Section~\ref{sec:prelim} contains the notation and auxiliary estimates.
Section~\ref{sec:universal} concerns general graphs.
The complete-graph and star-graph results are proved in Sections~\ref{sec:complete} and~\ref{sec:fracstar}, respectively.

\section{Preliminaries}\label{sec:prelim}
\subsection{Notation and basic reductions}
Throughout, $G=(V,E)$ is a finite connected simple undirected graph with $n=|V|\geq2$.
We denote its graph distance by $d_G$, its diameter by $D$, and the degree of $v$ by $\deg(v)$.
For integer $r\geq0$, let
\[
B_G(v,r)=\{w\in V:d_G(v,w)\leq r\}.
\]
Since $B_G(v,D)=V$, maxima over balls may be restricted to $0\leq r\leq D$.
We write $K_n$ for the complete graph. The star $S_n$ has center $0$, leaves $1,\ldots,n-1$, and edge set $\{\{0,i\}:1\leq i<n\}$.

\begin{dfn}\label{def:maximal}
For $0<\alpha<1$ and $f:V\to\R$, the centered fractional maximal operator is
\[
M_{\alpha,G}f(v)=\max_{0\leq r\leq D}
|B_G(v,r)|^{\alpha-1}\sum_{w\in B_G(v,r)}|f(w)|.
\]
\end{dfn}
Set $c=n^{\alpha-1}$ and $\sigma_f=\sum_{v\in V}f(v)$.
When vertices are indexed, we write $f_i=f(v_i)$.
\begin{dfn}\label{def:variation}
The edge variation and its associated sharp constants are defined by
\begin{align*}
\Var_p(f)&=\biggl(\sum_{\{u,v\}\in E}|f(u)-f(v)|^p\biggr)^{1/p}
&& (0<p<\infty),\\
\Var_\infty(f)&=\max_{\{u,v\}\in E}|f(u)-f(v)|,\\
\C_{\alpha,G;p,q}
&=\sup_{\Var_p(f)>0}\frac{\Var_q(M_{\alpha,G}f)}{\Var_p(f)},
&\C_{\alpha,G,p}&=\C_{\alpha,G;p,p}.
\end{align*}
Here $0<p\leq q\leq\infty$, and each undirected edge is counted once.
We use $1/\infty=0$ and $t_+=\max\{t,0\}$.
Connectedness implies that $\Var_p(f)=0$ exactly when $f$ is constant.
\end{dfn}

\begin{lem}[Nonnegative reduction and homogeneity]\label{lem:reduce}
The supremum defining $\C_{\alpha,G;p,q}$ may be restricted to nonnegative functions.
For $t>0$,
\[
M_{\alpha,G}(tf)=tM_{\alpha,G}f,\qquad \Var_p(tf)=t\Var_p(f).
\]
\end{lem}
Indeed, $M_{\alpha,G}f=M_{\alpha,G}|f|$ and $\Var_p(|f|)\leq\Var_p(f)$.
If $|f|\equiv a$, then $M_{\alpha,G}f\equiv n^\alpha a$, so its variation vanishes.

In particular, positive rescaling preserves every variation ratio.
Unlike the ordinary maximal operator, $M_{\alpha,G}$ does not commute with addition of constants; below, constant terms are retained in the ball sums.

\subsection{Auxiliary inequalities}
We record the standard finite-sequence inequalities used below, without proof.
\begin{lem}[Power and norm comparisons]\label{lem:powers}\label{lem:norm}
For $a_1,\ldots,a_k\geq0$,
\[
\biggl(\sum_{j=1}^k a_j\biggr)^p
\leq\sum_{j=1}^k a_j^p
\leq k^{1-p}\biggl(\sum_{j=1}^k a_j\biggr)^p
\qquad(0<p\leq1).
\]
For $1\leq p<\infty$, both inequalities are reversed.
For $\mathbf z\in\R^N$ and $0<p\leq q\leq\infty$,
\[
\|\mathbf z\|_q\leq\|\mathbf z\|_p
\leq N^{1/p-1/q}\|\mathbf z\|_q.
\]
If $0<p<q<\infty$, then
\[
\|\mathbf z\|_q\leq
\|\mathbf z\|_\infty^{1-p/q}\|\mathbf z\|_p^{p/q}.
\]
\end{lem}
The notation $\|\cdot\|_p$ denotes the usual unnormalized $\ell^p$ norm or quasi-norm.

\begin{dfn}\label{def:osc}
Write
\[
\osc(f)=\max_{v\in V}f(v)-\min_{v\in V}f(v),\qquad
\partial_EA=\{\{u,v\}\in E:u\in A,\ v\notin A\}
\]
for the oscillation and the edge boundary, respectively.
\end{dfn}
Thus $\Var_\infty(f)\leq\osc(f)$, with equality on $K_n$.

\begin{lem}[Oscillation and ordered gaps]\label{lem:oscinput}\label{lem:ordered-gaps}
For every real-valued $f$ and $0<p\leq\infty$,
\[
\osc(f)\leq D^{(1-1/p)_+}\Var_p(f).
\]
If $t_1\leq\cdots\leq t_n$ are the ordered vertex values and $\delta_j=t_{j+1}-t_j$, then
\[
\sum_{j=1}^{n-1}\delta_j^p\leq\Var_p(f)^p
\quad(1\leq p<\infty),\qquad
\max_j\delta_j\leq\Var_\infty(f).
\]
On $K_n$, one also has
\[
\Var_p(f)^p\geq(n-1)\osc(f)^p\qquad(0<p\leq1).
\]
\end{lem}
\begin{proof}
The first estimate follows by summing differences along a shortest path between a minimizing and a maximizing vertex, and applying the power inequalities or H\"older's inequality.
For the gap estimates, label the vertices so that $f(v_i)=t_i$ and put $A_j=\{v_1,\ldots,v_j\}$.
Connectedness gives $|\partial_EA_j|\geq1$. For $p\geq1$,
\[
\Var_p(f)^p
=\sum_{\substack{a<b\\\{v_a,v_b\}\in E}}
\biggl(\sum_{j=a}^{b-1}\delta_j\biggr)^p
\geq\sum_{j=1}^{n-1}|\partial_EA_j|\delta_j^p
\geq\sum_{j=1}^{n-1}\delta_j^p.
\]
An edge crossing $A_j$ has difference at least $\delta_j$, proving the infinite-exponent estimate.
For the last assertion, choose vertices $v_-,v_+$ attaining the minimum and maximum.
The edge between them and the two edges joining each remaining vertex to them are distinct.
Subadditivity of $t^p$ therefore gives
\[
\Var_p(f)^p\geq \osc(f)^p
+\sum_{v\notin\{v_-,v_+\}}
\bigl[(f(v)-\min f)^p+(\max f-f(v))^p\bigr]
\geq(n-1)\osc(f)^p.
\]
\end{proof}

\begin{lem}[Extreme-point principle]\label{lem:vertex}
A convex function on a nonempty compact convex polytope attains a maximum at an extreme point.
\end{lem}

\section{Two-sided estimates for individual finite graphs}\label{sec:universal}

We combine oscillation estimates, point-mass tests, and the relation between edge count and diameter to obtain two-sided bounds for a fixed graph.
\subsection{Oscillation comparisons and their parameter dependence}\label{sec:osc}
We compare a maximizing ball with the whole-graph ball.
Retaining the number of vertices in that ball gives an oscillation bound, which can then be converted into a variation bound using the tools of Section 2.
For $0<\alpha<1$, put
\[
c=n^{\alpha-1},\qquad
A_\alpha(n)=\max_{1\leq k<n}\bigl(k^\alpha-kc\bigr).
\]
The maximum is taken over integers $k=1,\ldots,n-1$.

\begin{lem}[Oscillation comparison for the maximal function]\label{thm:osc}
Let $n\geq2$ and $0<\alpha<1$. For every connected graph on $n$ vertices and every $f\geq0$,
\[
\osc(M_{\alpha,G}f)\leq A_\alpha(n)\osc(f).
\]
\end{lem}
\begin{proof}
Let $B$ be a ball attaining $\max_v M_{\alpha,G}f(v)$, and put $k=|B|$.
The whole-graph ball gives $M_{\alpha,G}f(v)\geq c\sum_w f(w)$ at every vertex.
If $k=n$, the maximal function is constant. Otherwise,
\begin{align*}
\osc(M_{\alpha,G}f)
&\leq k^{\alpha-1}\sum_{w\in B}f(w)-c\sum_{w\in V}f(w)\\
&=(k^\alpha-n^\alpha)\min f
 +(k^{\alpha-1}-c)\sum_{w\in B}(f(w)-\min f)\\
&\quad-c\sum_{w\notin B}(f(w)-\min f)\\
&\leq(k^\alpha-kc)\osc(f)\leq A_\alpha(n)\osc(f).
\end{align*}
Here the first and third terms in the decomposition are nonpositive, since $f\geq0$ and $k<n$.
\end{proof}

\begin{thm}[Optimal uniform oscillation coefficient for the fractional maximal function]\label{thm:sharp-osc-variation}
Fix $n\geq2$ and $0<\alpha<1$, and put $c=n^{\alpha-1}$. For integers $1\leq k,j<n$, define
\[
w_{k,j}=\begin{cases}
k^\alpha-(n-j)c,&j\leq n-k,\\
(n-j)(k^{\alpha-1}-c),&j>n-k.
\end{cases}
\]
Using the positive-part notation, define
\[
\Omega_{\alpha,p}(n)=\begin{cases}
A_\alpha(n),&0<p\leq1,\\
\displaystyle\max_{1\leq k<n}\biggl(\sum_{j=1}^{n-1}(w_{k,j})_+^{p'}\biggr)^{1/p'},
&1<p<\infty,\quad p'=p/(p-1),\\
\displaystyle\max_{1\leq k<n}\sum_{j=1}^{n-1}(w_{k,j})_+,&p=\infty.
\end{cases}
\]
Then, for every connected simple graph $G$ on $n$ vertices, every real-valued function $f$, and every $0<p\leq\infty$,
\[
\osc(M_{\alpha,G}f)\leq\Omega_{\alpha,p}(n)\Var_p(f).
\]
Moreover, this coefficient is optimal over the class $\mathscr G_n$ of all connected simple graphs on $n$ vertices:
\[
\sup_{G\in\mathscr G_n}\ \sup_{\Var_p(f)>0}
\frac{\osc(M_{\alpha,G}f)}{\Var_p(f)}=\Omega_{\alpha,p}(n).
\]
Equality is attained by functions on a path.
\end{thm}
\begin{proof}
By Lemma~\ref{lem:reduce}, it suffices to prove the upper bound for $f\geq0$.
Order its vertex values as $t_1\leq\cdots\leq t_n$ and put $\delta_j=t_{j+1}-t_j$.
The identities
\[
\sum_{i=1}^n t_i=nt_1+\sum_{j=1}^{n-1}(n-j)\delta_j,\qquad
\sum_{i=n-k+1}^n t_i=kt_1+\sum_{j=1}^{n-1}\min\{k,n-j\}\delta_j
\]
give the required coefficients.
Indeed, choose a ball attaining the maximum of $M_{\alpha,G}f$.
If this ball is $V$, the maximal function is constant.
Otherwise, writing $k<n$ for its cardinality and comparing its sum with the largest $k$ vertex values, we obtain
\begin{align*}
\osc(M_{\alpha,G}f)
&\leq k^{\alpha-1}\sum_{i=n-k+1}^n t_i-c\sum_{i=1}^n t_i\\
&=(k^\alpha-n^\alpha)t_1+\sum_{j=1}^{n-1}w_{k,j}\delta_j
\leq\sum_{j=1}^{n-1}(w_{k,j})_+\delta_j.
\end{align*}
For $1<p<\infty$, H\"older's inequality and Lemma~\ref{lem:ordered-gaps} imply
\[
\sum_{j=1}^{n-1}(w_{k,j})_+\delta_j
\leq\biggl(\sum_{j=1}^{n-1}(w_{k,j})_+^{p'}\biggr)^{1/p'}
\biggl(\sum_{j=1}^{n-1}\delta_j^p\biggr)^{1/p}
\leq\Omega_{\alpha,p}(n)\Var_p(f).
\]
For $p=\infty$, use $\delta_j\leq\Var_\infty(f)$.
For $0<p\leq1$, the weights increase up to $j=n-k$ and decrease afterwards, so
\[
(w_{k,j})_+\leq w_{k,n-k}=k^\alpha-kc\leq A_\alpha(n).
\]
Together with $\sum_j\delta_j=\osc(f)\leq\Var_p(f)$, this proves the upper bound in every case.

To prove sharpness, take the path with edges $\{v_j,v_{j+1}\}$, $1\leq j<n$, and define
\[
f(v_1)=0,\qquad f(v_i)=\sum_{j=1}^{i-1}\delta_j\quad(2\leq i\leq n)
\]
for nonnegative gaps that are not all zero.
Then $\Var_p(f)=\|\boldsymbol\delta\|_p$.
The balls centered at $v_1$ are initial segments.
Their ordinary averages are nondecreasing, and multiplication by their cardinality to the power $\alpha$ preserves this property.
Thus
\[
M_{\alpha,G}f(v_1)=c\sum_{i=1}^n f(v_i)
=\min_v M_{\alpha,G}f(v).
\]
The ball of radius $k-1$ centered at $v_n$ consists of the last $k$ vertices. Consequently,
\[
\osc(M_{\alpha,G}f)\geq
k^{\alpha-1}\sum_{i=n-k+1}^n f(v_i)-c\sum_{i=1}^n f(v_i)
=\sum_{j=1}^{n-1}w_{k,j}\delta_j.
\]
If $1<p<\infty$, choose $k$ attaining the maximum defining $\Omega_{\alpha,p}(n)$ and set
$\delta_j=(w_{k,j})_+^{1/(p-1)}$.
Since $w_{k,n-k}>0$, the denominator below is positive, and
\[
\frac{\osc(M_{\alpha,G}f)}{\Var_p(f)}
\geq\frac{\sum_j(w_{k,j})_+^{p'}}
{\bigl(\sum_j(w_{k,j})_+^{p'}\bigr)^{1/p}}
=\Omega_{\alpha,p}(n).
\]
For $p=\infty$, take $\delta_j=\1_{\{w_{k,j}>0\}}$ for a maximizing $k$.
Then $\Var_\infty(f)=1$ and the same lower bound is $\sum_j(w_{k,j})_+$.
For $0<p\leq1$, choose $k$ with $k^\alpha-kc=A_\alpha(n)$ and take
$\delta_{n-k}=1$, with all other gaps zero.
The original function has variation $1$, and the preceding lower bound is $w_{k,n-k}=A_\alpha(n)$.
In each case, the upper bound already proved forces equality.
\end{proof}
\begin{rem}\label{rem:refined-osc-factor}
Since $A_\alpha(n)\geq1-c>0$, define
\[
R_{\alpha,p}(n,D)=\min\biggl\{D^{(1-1/p)_+},\frac{\Omega_{\alpha,p}(n)}{A_\alpha(n)}\biggr\}.
\]
Combining Lemmas~\ref{thm:osc} and~\ref{lem:oscinput} with Theorem~\ref{thm:sharp-osc-variation} gives
\[
\osc(M_{\alpha,G}f)\leq A_\alpha(n)R_{\alpha,p}(n,D)\Var_p(f).
\]
For $p\leq1$, one has $R_{\alpha,p}(n,D)=1$.
\end{rem}

\subsection{Point-mass tests, diameter counts, and edge-difference estimates}
\begin{lem}[Improved lower bounds from point masses]\label{lem:simple-dirac-lower}
Let $n\geq4$ and $0<\alpha<1$. For each $0<p\leq\infty$, there is a vertex $w$ such that
\[
\frac{\Var_p(M_{\alpha,G}\delta_w)}{\Var_p(\delta_w)}
\geq 1-3^{\alpha-1}.
\]
If $0<p\leq1$, the vertex can be chosen so that
\begin{equation}
\frac{\Var_p(M_{\alpha,G}\delta_w)}{\Var_p(\delta_w)}
\geq 1-4^{\alpha-1}. \label{eq:l1}
\end{equation}
Here $\delta_w(w)=1$, and $\delta_w$ vanishes at all other vertices.
\end{lem}
\begin{proof}
For any point mass, nestedness of balls and $\alpha-1<0$ give
\[
M_{\alpha,G}\delta_w(w)=1,\qquad
M_{\alpha,G}\delta_w(v)=|B_G(v,d_G(v,w))|^{\alpha-1}\quad(v\ne w).
\]
In particular, if $v$ is adjacent to $w$,
\[
M_{\alpha,G}\delta_w(v)=(\deg(v)+1)^{\alpha-1}.
\]
The function $\delta_w$ has absolute edge difference $1$ precisely on the $\deg(w)$ edges incident to $w$. Thus
\[
\Var_p(\delta_w)=\deg(w)^{1/p}\quad(0<p<\infty),
\qquad \Var_\infty(\delta_w)=1.
\]
\paragraph{The lower bound for all exponents.}
Choose a vertex $w$ of minimum degree, and put $\delta(G)=\deg(w)$.
If $\delta(G)\geq2$, every neighbor $v$ satisfies $\deg(v)\geq \delta(G)$, so
\[
|M_{\alpha,G}\delta_w(w)-M_{\alpha,G}\delta_w(v)|
=1-(\deg(v)+1)^{\alpha-1}
\geq1-(\delta(G)+1)^{\alpha-1}
\geq1-3^{\alpha-1}.
\]
For finite $p$, summation over these $\delta(G)$ edges gives
\[
\Var_p(M_{\alpha,G}\delta_w)^p\geq \delta(G)(1-3^{\alpha-1})^p.
\]
Divide by $\Var_p(\delta_w)^p=\delta(G)$ and take the $1/p$ power to obtain the first inequality.
For the infinite exponent, any one of these edges suffices.
If $\delta(G)=1$, let $v$ be the unique neighbor. If $\deg(v)=1$, then $\{w,v\}$ has no edge to the remaining vertices, contradicting connectedness and $n\geq4$.
Thus $\deg(v)\geq2$, and the same edge gives $|M_{\alpha,G}\delta_w(w)-M_{\alpha,G}\delta_w(v)|\geq1-3^{\alpha-1}$, while the original variation is $1$ for every exponent.

\paragraph{The lower bound for small exponents.}
Assume $0<p\leq1$.
Let $\mathcal P$ be a simple path from $w$ to a vertex $z$ with $M_{\alpha,G}\delta_w(z)\leq4^{\alpha-1}$, and denote its edge set by $\mathcal E(\mathcal P)$. Lemma~\ref{lem:powers} gives
\begin{align*}
(1-4^{\alpha-1})^p
&\leq |M_{\alpha,G}\delta_w(w)-M_{\alpha,G}\delta_w(z)|^p\\
&\leq\biggl(\sum_{\{u,v\}\in\mathcal E(\mathcal P)}|M_{\alpha,G}\delta_w(u)-M_{\alpha,G}\delta_w(v)|\biggr)^p\\
&\leq\sum_{\{u,v\}\in\mathcal E(\mathcal P)}|M_{\alpha,G}\delta_w(u)-M_{\alpha,G}\delta_w(v)|^p.
\end{align*}
Consequently, if there are $\delta(G)$ such paths with pairwise disjoint edge sets, then
\[
\Var_p(M_{\alpha,G}\delta_w)^p\geq\sum_{j=1}^{\delta(G)}
\sum_{\{u,v\}\in\mathcal E(\mathcal P_j)}|M_{\alpha,G}\delta_w(u)-M_{\alpha,G}\delta_w(v)|^p
\geq \delta(G) (1-4^{\alpha-1})^p.
\]
We construct the required edges or paths according to the minimum degree.

If the minimum degree $\delta(G)\geq3$, choose any vertex $w$ of minimum degree. Each of its $\delta(G)$ neighbors satisfies $M_{\alpha,G}\delta_w(v)\leq(\delta(G)+1)^{\alpha-1}\leq4^{\alpha-1}$, so the $\delta(G)$ edges incident to $w$ are the required paths.

If the minimum degree is $1$, choose a leaf $w$ and denote its neighbor by $v$.
When $\deg(v)\geq3$, the single-edge path $(w,v)$ suffices.
When $\deg(v)=2$, denote its other neighbor by $u$.
The vertex $u$ must have a neighbor $z\ne v$; otherwise $\{w,v,u\}$ has no edge to the remaining vertices.
Since $w$ is a leaf, $z\ne w$, so $w,v,u,z$ are distinct.
Also, $d_G(u,w)=2$, and
\[
\{w,v,u,z\}\subset B_G(u,2),\qquad
M_{\alpha,G}\delta_w(u)=|B_G(u,2)|^{\alpha-1}\leq4^{\alpha-1}.
\]
Thus the path $(w,v,u)$ has the required property.

It remains to consider minimum degree $2$. If there is a vertex of degree at least $3$, choose a path from a degree-$2$ vertex to it.
Denote the first edge on this path leaving $\{v:\deg(v)=2\}$ by $\{w,v\}$, where
$\deg(w)=2$, $\deg(v)\geq3$.
Let $u$ be the other neighbor of $w$. We already have $M_{\alpha,G}\delta_w(v)\leq4^{\alpha-1}$.
If $\deg(u)\geq3$, use the two single-edge paths $(w,v)$ and $(w,u)$.
If $\deg(u)=2$, let $z$ be its other neighbor.
When $z=v$, take
\[
\mathcal P_1=(w,v),\qquad \mathcal P_2=(w,u,v).
\]
Their edge sets are $\{\{w,v\}\}$ and $\{\{w,u\},\{u,v\}\}$, respectively, and are disjoint.
When $z\ne v$, the vertices $z,w$ are not adjacent, since the only neighbors of $w$ are $u,v$. Hence $d_G(z,w)=2$.
Choose a neighbor $t\ne u$ of $z$, which exists because the minimum degree is $2$. The preceding nonadjacency gives $t\ne w$.
Thus $w,u,z,t$ are distinct, and
\[
\{w,u,z,t\}\subset B_G(z,2),\qquad M_{\alpha,G}\delta_w(z)\leq4^{\alpha-1}.
\]
Take $\mathcal P_1=(w,v)$ and $\mathcal P_2=(w,u,z)$; their edge sets are again disjoint.

If every vertex has degree $2$, then $G$ is a cycle.

Label its vertices cyclically as $v_0,\ldots,v_{n-1}$.
Put $w=v_0$ and $z=v_2$. Then $d_G(w,z)=2$, and $B_G(z,2)$ contains at least $v_0,v_1,v_2,v_3$, so $M_{\alpha,G}\delta_w(z)\leq4^{\alpha-1}$.
Take the paths in the two directions around the cycle:
\[
\mathcal P_1=(v_0,v_1,v_2),\qquad
\mathcal P_2=(v_0,v_{n-1},v_{n-2},\ldots,v_2).
\]
They are edge-disjoint and both satisfy the required endpoint condition.

In every case, the chosen $w$ has minimum degree $\delta(G)$, and $\Var_p(M_{\alpha,G}\delta_w)^p\geq \delta(G) (1-4^{\alpha-1})^p$.
Together with $\Var_p(\delta_w)^p=\delta(G)$, this proves the second inequality.

\end{proof}

\begin{lem}[Diameter, edge count, and a point-mass test]\label{lem:diam-simple}
Let $D=\diam(G)$, and define the quantity depending only on $n,D$ by
\[
N_E(n,D)=\min\biggl\{\frac{n(n-1)}2,
D+\frac{(n-D-1)(n-D+4)}2\biggr\}.
\]
Then $|E|\leq N_E(n,D)$.
For $0<\alpha<1$ and $0<p\leq q\leq\infty$, there is a vertex $w$ such that $f=\delta_w$ satisfies
\begin{equation}
\frac{\Var_q(M_{\alpha,G}f)}{\Var_p(f)}
\geq\frac{1-n^{\alpha-1}}
{(n-D)^{1/p}D^{\max(1-1/q,0)}}. \label{eq:l2}
\end{equation}
\end{lem}
\begin{proof}
The edge-count estimate is Ore's bound \cite{Ore}.
For the variation estimate, choose a shortest path $v_0,\ldots,v_D$ with $d_G(v_0,v_D)=D$, and take $f=\delta_{v_0}$.
No vertex $v_2,\ldots,v_D$ is adjacent to $v_0$, so
\[
\Var_p(f)=\deg(v_0)^{1/p}\leq(n-D)^{1/p}.
\]
The point-mass formula gives
$M_{\alpha,G}f(v_0)=1$ and $M_{\alpha,G}f(v_D)=n^{\alpha-1}$.
Summing differences along the chosen path and applying the finite-sequence inequalities yields
\[
1-n^{\alpha-1}\leq
D^{(1-1/q)_+}\Var_q(M_{\alpha,G}f).
\]
These two estimates prove \eqref{eq:l2}, with the usual convention for infinite exponents.
\end{proof}

We next combine the edge-count estimate with convexity to control the sum of powered edge differences at a prescribed oscillation.
\begin{lem}[Powered edge-difference sums controlled by oscillation]\label{lem:energy-size}
Let $n\geq2$, $D=\diam(G)$, and $0<q\leq\infty$. Define
\[
Q_q(n,D)=
\begin{cases}
N_E(n,D)^{1-q}\min\{N_E(n,D),\lfloor n^2/4\rfloor\}^{q},&0<q<1,\\
\min\{N_E(n,D),\lfloor n^2/4\rfloor\},&1\leq q<\infty,\\
1,&q=\infty.
\end{cases}
\]
Every function $h$ satisfies $\Var_q(h)\leq Q_q(n,D)^{1/q}\osc(h)$.
At the infinite endpoint, the factor is interpreted as $1$.
\end{lem}
\begin{proof}
The case $\osc(h)=0$ is immediate.
Otherwise, normalize the vertex values by $g_i=(h(v_i)-\min h)/\osc(h)$.
For $1\leq q<\infty$, the convex function
\[
\Phi_q(\mathbf z)=\sum_{i<j}|z_i-z_j|^q
\]
attains its maximum on $[0,1]^n$ at a vertex of the cube.
At such a vertex, with $k$ coordinates equal to $1$, its value is $k(n-k)\leq\lfloor n^2/4\rfloor$.
Since each edge difference is also at most $\osc(h)$, it follows that
\[
\Var_q(h)^q\leq
\min\{|E|,\lfloor n^2/4\rfloor\}\osc(h)^q
\leq Q_q(n,D)\osc(h)^q
\qquad(1\leq q<\infty).
\]
For $0<q<1$, the power inequality and the case $q=1$ give
\[
\Var_q(h)^q\leq |E|^{1-q}\Var_1(h)^q
\leq N_E(n,D)^{1-q}Q_1(n,D)^q\osc(h)^q.
\]
The infinite-exponent estimate is $\Var_\infty(h)\leq\osc(h)$.
\end{proof}

\subsection{Equal-exponent and mixed estimates for the fractional operator}\label{sec:fracgeneral}
We prove two-sided estimates for an individual graph, first for equal exponents and then for mixed exponents.
Fix $0<\alpha<1$ and put $c=n^{\alpha-1}$. The quantities $A_\alpha(n)$ and $Q_r(n,D)$ were defined in Subsection~\ref{sec:osc} and Lemma~\ref{lem:energy-size}, respectively.

\subsubsection{Two-sided bounds for equal exponents}
Combining the oscillation upper bound with the point-mass lower bounds gives the following estimate.
\begin{thm}[Two-sided fractional bounds for equal exponents]\label{thm:general-pp}
Let $n\geq4$ and $0<\alpha<1$. Define
\[
b_{\alpha,p}=\begin{cases}1-4^{\alpha-1},&0<p\leq1,\\1-3^{\alpha-1},&1<p\leq\infty.\end{cases}
\]
Then
\begin{equation}
\ell_{\alpha,p}(n,D)
\leq\C_{\alpha,G,p}
\leq A_\alpha(n)Q_p(n,D)^{1/p}R_{\alpha,p}(n,D), \label{eq:fgpp}
\end{equation}
where
\[
\ell_{\alpha,p}(n,D)=\begin{cases}
\max\{b_{\alpha,p},(1-c)(n-D)^{-1/p}D^{-(1-1/p)_+}\},&p<\infty,\\
\max\{b_{\alpha,\infty},1-n^{-(1-\alpha)/D}\},&p=\infty.
\end{cases}
\]

\end{thm}
\begin{proof}
For finite $p$, the lower bound follows from Lemmas~\ref{lem:simple-dirac-lower} and~\ref{lem:diam-simple}.
For $p=\infty$, put $\Delta=\max_v\deg(v)$.
Choose a vertex $v$ of degree $\Delta$ and a neighbor $w$.
The test $f=\delta_w$ gives
\[
\C_{\alpha,G,\infty}\geq1-(\Delta+1)^{\alpha-1}.
\]
The elementary ball-growth estimate
$|B_G(u,r+1)|\leq(\Delta+1)|B_G(u,r)|$
implies $n\leq(\Delta+1)^D$. Hence
\[
\C_{\alpha,G,\infty}\geq1-n^{-(1-\alpha)/D}.
\]
Together with Lemma~\ref{lem:simple-dirac-lower}, this proves the lower bound.
For the upper bound, Lemma~\ref{lem:energy-size} and Remark~\ref{rem:refined-osc-factor} give
\[
\Var_p(M_{\alpha,G}f)
\leq Q_p(n,D)^{1/p}\osc(M_{\alpha,G}f)
\leq A_\alpha(n)Q_p(n,D)^{1/p}R_{\alpha,p}(n,D)\Var_p(f),
\]
where the counting factor equals $1$ when $p=\infty$.
\end{proof}

\subsubsection{Two-sided bounds for mixed exponents}
The degree factor of a point mass refines the mixed-exponent lower bound. The upper bound again follows from oscillation comparison.
\begin{thm}[Two-sided fractional bounds for mixed exponents]\label{thm:generalbounds}
Let $n\geq4$, $0<\alpha<1$, and $0<p<q\leq\infty$. Then
\begin{equation}
\max\biggl\{L_{\alpha,p,q}(n),
\frac{1-c}{(n-D)^{1/p}D^{(1-1/q)_+}}\biggr\}
\leq\C_{\alpha,G;p,q}
\leq A_\alpha(n)Q_q(n,D)^{1/q}R_{\alpha,p}(n,D), \label{eq:fgmix}
\end{equation}
where
\[
L_{\alpha,p,q}(n)=
\begin{cases}
\min\{3^{1/q-1/p}(1-4^{\alpha-1}),(n-1)^{1/q-1/p}(1-c)\},&q\leq1,\\
\min\{2^{1/q-1/p}(1-3^{\alpha-1}),(n-1)^{1/q-1/p}(1-c)\},&q>1.
\end{cases}
\]
Bounds with the diameter eliminated are given in Corollary~\ref{cor:general-geometryfree}.
\end{thm}
\begin{proof}
For the common lower bound, we compare point-mass variation ratios according to the minimum degree; the upper bound again follows from oscillation. Since $p<q\leq\infty$, the exponent $p$ is finite and positive, and $1/q-1/p<0$. We use $1/q=0$ throughout when $q=\infty$.

\paragraph{Endpoint minima of an auxiliary function.}
Define
\[
F(t)=t^{1/q-1/p}\bigl[1-(t+1)^{\alpha-1}\bigr],\qquad t>0.
\]
Since $F(t)>0$, logarithmic differentiation gives
\[
\frac{tF'(t)}{F(t)}=1/q-1/p+H(t),\qquad
H(t)=\frac{(1-\alpha)t}{(t+1)[(t+1)^{1-\alpha}-1]}.
\]
Differentiating $H$ gives
\[
H'(t)=\frac{(1-\alpha)\bigl[(t+1)^{1-\alpha}-1
-(1-\alpha)t(t+1)^{1-\alpha}\bigr]}
{(t+1)^2\bigl[(t+1)^{1-\alpha}-1\bigr]^2}.
\]
The denominator is positive. Since $u^{-\alpha}\leq1$ for $u\geq1$,
\[
(t+1)^{1-\alpha}-1
=\int_1^{t+1}(1-\alpha)u^{-\alpha}\,du
\leq(1-\alpha)t.
\]
The bracket in the numerator is at most
\[
(1-\alpha)t\bigl[1-(t+1)^{1-\alpha}\bigr]<0.
\]
Thus $H$ is strictly decreasing, and $F'$ changes sign at most once, from positive to negative.
It follows that
\[
F(t)\geq\min\{F(a),F(b)\}\qquad(0<a\leq t\leq b).
\]

\paragraph{Degree-dependent lower bounds.}
Denote the minimum degree by $\delta(G)$. First suppose $q>1$, including $q=\infty$.
If $\delta(G)\geq2$, choose a minimum-degree vertex $w$. The point-mass calculation in Lemma~\ref{lem:simple-dirac-lower} shows that each neighbor $v$ satisfies
\[
|M_{\alpha,G}\delta_w(w)-M_{\alpha,G}\delta_w(v)|
\geq1-(\delta(G)+1)^{\alpha-1}.
\]
For finite $q$, sum over these $\delta(G)$ edges; for the infinite exponent, use any one of them. Both cases give
\[
\frac{\Var_q(M_{\alpha,G}\delta_w)}{\Var_p(\delta_w)}
\geq \delta(G)^{1/q-1/p}\bigl[1-(\delta(G)+1)^{\alpha-1}\bigr]
=F(\delta(G)).
\]
Applying the endpoint estimate on $2\leq \delta(G)\leq n-1$ gives
\[
F(\delta(G))\geq\min\{F(2),F(n-1)\}
=\min\{2^{1/q-1/p}(1-3^{\alpha-1}),(n-1)^{1/q-1/p}(1-c)\}.
\]
If $\delta(G)=1$, take the indicator of a leaf. Its neighbor has degree at least $2$, and the original $p$-variation is $1$, so the ratio is at least $1-3^{\alpha-1}$.
Since $1/q-1/p<0$, we have $1-3^{\alpha-1}\geq2^{1/q-1/p}(1-3^{\alpha-1})$, and the same bound follows.

Now suppose $q\leq1$, so $p\leq q\leq1$.
If $\delta(G)\geq3$, the same point-mass calculation bounds the ratio below by $F(\delta(G))$. Hence
\[
F(\delta(G))\geq\min\{F(3),F(n-1)\}
=\min\{3^{1/q-1/p}(1-4^{\alpha-1}),(n-1)^{1/q-1/p}(1-c)\}.
\]
If $\delta(G)\in\{1,2\}$, choose the minimum-degree vertex $w$ constructed in the proof of Lemma~\ref{lem:simple-dirac-lower}.
That construction gives $\delta(G)$ edge-disjoint paths starting at $w$, with maximal-function values at their terminal vertices no greater than $4^{\alpha-1}$.
Taking the power exponent in that proof to be $q$ gives
\[
\Var_q(M_{\alpha,G}\delta_w)^q\geq \delta(G)(1-4^{\alpha-1})^q,
\qquad \Var_p(\delta_w)=\delta(G)^{1/p}.
\]
Therefore
\[
\frac{\Var_q(M_{\alpha,G}\delta_w)}{\Var_p(\delta_w)}
\geq \delta(G)^{1/q-1/p}(1-4^{\alpha-1})
\geq3^{1/q-1/p}(1-4^{\alpha-1}).
\]
The final step uses $\delta(G)\leq2<3$ and $1/q-1/p<0$.
This proves the common lower bound $L_{\alpha,p,q}(n)$ in each exponent range.
Lemma~\ref{lem:diam-simple} also gives the diameter lower bound, and the maximum of the two may be taken.

\paragraph{The upper bound.}
By Lemma~\ref{lem:reduce}, it suffices to consider $f\geq0$.
Applying Lemma~\ref{lem:energy-size} and Remark~\ref{rem:refined-osc-factor}, we obtain
\begin{align*}
\Var_q(M_{\alpha,G}f)
&\leq Q_q(n,D)^{1/q}\osc(M_{\alpha,G}f)\\
&\leq Q_q(n,D)^{1/q}A_\alpha(n)R_{\alpha,p}(n,D)\Var_p(f).
\end{align*}
The counting factor is $1$ when $q=\infty$.
\end{proof}

\subsection{Bounds independent of geometric graph parameters}\label{sec:geometryfree}
Maximizing over all possible diameters gives an upper bound depending only on the number of vertices, the fractional parameter, and the variation exponents.
\begin{cor}[Parameter-only two-sided bounds]\label{cor:general-geometryfree}
Let $n\geq4$. Put $\mathcal H_{\alpha;p,q}(n)=\max_{1\leq r<n}Q_q(n,r)^{1/q}R_{\alpha,p}(n,r)$, where the integer $r$ ranges over $1,\ldots,n-1$.
Then
\[
b_{\alpha,p}\leq\C_{\alpha,G,p}\leq A_\alpha(n)\mathcal H_{\alpha;p,p}(n),
\]
\[
L_{\alpha,p,q}(n)\leq\C_{\alpha,G;p,q}\leq A_\alpha(n)\mathcal H_{\alpha;p,q}(n)\qquad(p<q).
\]
\end{cor}
\begin{proof}
Apply Theorems~\ref{thm:general-pp} and~\ref{thm:generalbounds}, using $1\leq D\leq n-1$.
\end{proof}

\section{Sharp constants on complete graphs}\label{sec:complete}

Fix $n\geq2$, label the vertices of $K_n$ by $v_1,\ldots,v_n$, and write
$f_i=f(v_i)$, $\sigma_f=\sum_{i=1}^n f_i$.
The only balls in a complete graph are singletons and the whole vertex set. Thus, for $f\geq0$,
\[
M_{\alpha,K_n}f(v_i)=\max\{f_i,c\sigma_f\},
\qquad c=n^{\alpha-1},\qquad \frac1n<c<1.
\]
We use the known mean-truncation inequality and then raise the truncation level from $1/n$ to $c$.

\subsection{The mean-truncation inequality}
The following inequality is a known result for the ordinary maximal operator on complete graphs; see \cite{GKM}*{Theorem 1}.
\begin{lem}[Normalized mean truncation]\label{lem:mean-truncation}
Let $x_i\geq0$, $\sum_{i=1}^n x_i=1$, and $n\geq2$. For $0<p<\infty$, we have
\[
\sum_{i<j}|(x_j-1/n)_+-(x_i-1/n)_+|^p
\leq\biggl(1-\frac1n\biggr)^p\sum_{i<j}|x_j-x_i|^p.
\]
\end{lem}
\begin{proof}
For $n\geq3$, apply \cite{GKM}*{Theorem 1} to a nonnegative function of total sum $1$ on $K_n$.
For $n=2$, the assertion follows from
$|(x_2-1/2)_+-(x_1-1/2)_+|=|x_2-x_1|/2$.
\end{proof}

\subsection{Monotonicity at higher truncation levels}\label{sec:fraccomplete}
Raising the mean-truncation level to the fractional truncation level $c$ also requires control of the factor $1-c$.
This accounts for the denominator $(1-t)^p$ in the following lemma.

\begin{lem}[Monotonicity in the truncation level]\label{lem:threshold}
Let $n\geq2$, $x_1,\ldots,x_n\geq0$, $\sum_{i=1}^n x_i=1$, and $0<p<\infty$. Then
\[
F_p(t)=\frac{\sum_{i<j}|(x_j-t)_+-(x_i-t)_+|^p}{(1-t)^p},
\qquad 1/n\leq t<1
\]
is nonincreasing.
\end{lem}
\begin{proof}
Consider an open subinterval of $(1/n,1)$ containing none of the $x_i$.
On this interval, the set $\mathcal I=\{i:x_i>t\}$ and its cardinality $k\leq n-1$ are fixed.
If $k=0$, then $F_p=0$.
Otherwise, put $u_i=x_i-t>0$ for $i\in\mathcal I$ and
$\Phi=\sum_{i<j,\ i,j\in\mathcal I}|x_j-x_i|^p$.
The numerator of $F_p$ is $\Phi+(n-k)\sum_{i\in\mathcal I}u_i^p$, so
\begin{equation}\label{eq:k1}
F_p'(t)=\frac{p}{(1-t)^{p+1}}
\biggl[\Phi-(n-k)\sum_{i\in\mathcal I}(1-x_i)u_i^{p-1}\biggr].
\end{equation}
For every $i\in\mathcal I$, normalization gives
\[
1-x_i\geq(k-1)t,\qquad u_i\leq1-kt.
\]
Moreover, $|u_j-u_i|^p\leq u_i^p+u_j^p$ implies
$\Phi\leq(k-1)\sum_{i\in\mathcal I}u_i^p$.
Consequently, the expression in brackets in \eqref{eq:k1} is at most
\begin{align*}
&(k-1)\sum_{i\in\mathcal I}\bigl[u_i-(n-k)t\bigr]u_i^{p-1}\\
&\qquad\leq(k-1)(1-nt)\sum_{i\in\mathcal I}u_i^{p-1}\leq0.
\end{align*}
All powers $u_i^{p-1}$ are positive and finite on the interval under consideration, including when $p<1$.
Thus $F_p$ is nonincreasing between consecutive truncation levels.
Continuity across the finitely many levels $x_i$ proves the assertion on $[1/n,1)$.
\end{proof}

\subsection{Sharp constants for all exponent pairs}
For $n\geq4$, the case $p=1$ and $q=1/(1-\alpha)$ was obtained in \cite{LLW}*{Theorem 1.5}.
\begin{thm}[Sharp $(p,q)$ constants on complete graphs]\label{thm:complete-mixed}
Let $n\geq2$, $0<\alpha<1$, and $c=n^{\alpha-1}$.
For all $0<p\leq q\leq\infty$,
\[
\C_{\alpha,K_n;p,q}=(1-c)(n-1)^{1/q-1/p}.
\]
Equivalently, every real-valued function $f:V(K_n)\to\R$ satisfies
\[
\Var_q(M_{\alpha,K_n}f)
\leq(1-c)(n-1)^{1/q-1/p}\Var_p(f),
\]
and every point mass attains equality.
\end{thm}
\begin{proof}
By Lemma~\ref{lem:reduce}, it suffices to prove the upper bound for nonnegative $f$.
The maximal function of a constant function is constant, so assume that $f$ is nonconstant.
For $0<p<\infty$, set $g=f/\sigma_f$ and $x_i=g(v_i)$.
Since $M_{\alpha,K_n}g(v_i)=c+(x_i-c)_+$, Lemmas~\ref{lem:mean-truncation} and~\ref{lem:threshold} give
\[
\Var_p(M_{\alpha,K_n}g)^p
=(1-c)^pF_p(c)\leq(1-c)^pF_p(1/n)
\leq(1-c)^p\Var_p(g)^p.
\]
Homogeneity yields
\begin{equation}\label{eq:complete-diagonal-bound}
\Var_p(M_{\alpha,K_n}f)\leq(1-c)\Var_p(f).
\end{equation}

The whole-graph term is at least the minimum vertex value, and therefore
\begin{equation}\label{eq:complete-output-oscillation}
\Var_\infty(M_{\alpha,K_n}f)=(\max f-c\sigma_f)_+.
\end{equation}
Since $\sigma_f\geq \max f+(n-1)\min f$ and $nc>1$,
\[
\max f-c\sigma_f\leq(1-c)(\max f-\min f)+(1-nc)\min f\leq(1-c)(\max f-\min f).
\]
It follows that
\begin{equation}\label{eq:complete-infinity-diagonal}
\Var_\infty(M_{\alpha,K_n}f)\leq(1-c)\Var_\infty(f).
\end{equation}
For $0<p\leq1$, the complete-graph estimate in Lemma~\ref{lem:oscinput} gives
\[
\Var_\infty(M_{\alpha,K_n}f)\leq
(1-c)(n-1)^{-1/p}\Var_p(f).
\]
For $1<p<\infty$, choose $i_0$ with $f_{i_0}=\max f$.
The identity
\[
\frac{1-c}{n-1}\sum_{j\ne i_0}(\max f-f_j)-(\max f-c\sigma_f)
=\frac{nc-1}{n-1}(\sigma_f-\max f)\geq0
\]
and H\"older's inequality yield
\[
(\max f-c\sigma_f)_+
\leq\frac{1-c}{n-1}\sum_{j\ne i_0}(\max f-f_j)
\leq(1-c)(n-1)^{-1/p}\Var_p(f).
\]
Thus, for every $0<p\leq\infty$,
\begin{equation}\label{eq:complete-to-infinity}
\Var_\infty(M_{\alpha,K_n}f)
\leq(1-c)(n-1)^{-1/p}\Var_p(f).
\end{equation}

The cases $q=p$ and $q=\infty$ have been proved.
For $p<q<\infty$, interpolation of finite-sequence norms gives
\begin{align*}
\Var_q(M_{\alpha,K_n}f)
&\leq\Var_\infty(M_{\alpha,K_n}f)^{1-p/q}
       \Var_p(M_{\alpha,K_n}f)^{p/q}\\
&\leq(1-c)(n-1)^{1/q-1/p}\Var_p(f).
\end{align*}
Finally, if $f=\delta_{v_n}$, the maximal function equals $1$ at $v_n$ and $c$ elsewhere.
Thus, for every $0<r\leq\infty$,
\[
\Var_r(f)=(n-1)^{1/r},\qquad
\Var_r(M_{\alpha,K_n}f)=(1-c)(n-1)^{1/r}.
\]
Every point mass therefore attains the stated constant.
\end{proof}

\section{Variation estimates on star graphs}\label{sec:fracstar}
A leaf of a star graph has a two-point ball. Thus, besides the singleton value and the whole-graph term, its fractional maximal function contains a local combination of the center value and the leaf value.
We begin with the known sharp estimate for $p=1$ and then determine the exact parameter range of the $p=2$ formula. We then obtain the sharp infinite-exponent constant and use comparisons of variation exponents to deduce mixed estimates.

Throughout this section, let $n\geq3$ and $0<\alpha<1$. Denote the center of $S_n$ by $0$ and its leaves by $1,\ldots,n-1$, and write
\[
b=2^{\alpha-1},\qquad c=n^{\alpha-1}.
\]
Then $1/2<b<1$, $0<c<b$, and $nc=n^\alpha>2^\alpha=2b$.
By Lemma~\ref{lem:reduce}, in proving upper bounds we may assume
$f=(x,y_1,\ldots,y_{n-1})\geq0$. Set $\sigma_f=x+\sum_{i=1}^{n-1}y_i$.
The balls centered at the center are its singleton and the whole graph. Those centered at a leaf are its singleton, the two-point set consisting of that leaf and the center, and the whole graph. Hence
\begin{equation}\label{eq:f1}
\begin{split}
M_{\alpha,S_n}f(0)&=\max\{x,c\sigma_f\},\\
M_{\alpha,S_n}f(i)&=\max\{y_i,b(x+y_i),c\sigma_f\},\qquad 1\leq i\leq n-1.
\end{split}
\end{equation}

A single-leaf indicator gives the lower bound used repeatedly in this section. Take $f=\delta_1$, with $f(1)=1$ and zero at every other vertex.
By \eqref{eq:f1},
\[
M_{\alpha,S_n}\delta_1(1)=1,\qquad
M_{\alpha,S_n}\delta_1(i)=c\quad(i\ne1).
\]
The absolute edge-difference vectors of the original function and its maximal function are $(1,0,\ldots,0)$ and $(1-c,0,\ldots,0)$, respectively.
Thus, for all $0<p\leq q\leq\infty$,
\[
\Var_p(\delta_1)=1,\qquad
\Var_q(M_{\alpha,S_n}\delta_1)=1-c,\qquad
\C_{\alpha,S_n;p,q}\geq1-c.
\]

\subsection{The sharp constant for \texorpdfstring{$p=1$}{p=1}}
The following estimate is proved in \cite{LLW}*{proof of Theorem 1.6, (4.4)}.
We include a proof for completeness.
\begin{thm}[The sharp $p=1$ constant on star graphs]\label{thm:fracstar-p1-sharp}
For every $n\geq4$ and $0<\alpha<1$,
\begin{equation}\label{eq:fs4}
\C_{\alpha,S_n,1}=1-n^{\alpha-1}.
\end{equation}
Every single-leaf indicator attains equality.
\end{thm}
\begin{proof}
Only the upper bound remains to be proved. By the reduction to nonnegative functions and homogeneity, we may assume $f\geq0$ and $\sigma_f=1$; the zero function satisfies the conclusion separately.
Then $0\leq x\leq1$, $\sum_i y_i=1-x$, and
\[
M_{\alpha,S_n}f(0)=\max\{x,c\},\qquad
M_{\alpha,S_n}f(i)=\max\{y_i,b(x+y_i),c\}.
\]
Since $n\geq4$ and $0<\alpha<1$, we also have
\begin{equation}\label{eq:srel}
c\leq4^{\alpha-1}=b^2,\qquad nc\geq4^\alpha=4b^2.
\end{equation}
We prove that
\begin{equation}\label{eq:sgoal}
\Var_1(M_{\alpha,S_n}f)\leq(1-c)\Var_1(f).
\end{equation}

\paragraph{Case 1: $0\leq x\leq c$.}
The maximal function equals $c$ at the center and is at least $c$ at every leaf.
For fixed $x$, define
\[
\phi_x(y)=\max\{y-c,b(x+y)-c,0\}-(1-c)|y-x|.
\]
The desired inequality is equivalent to $\sum_i\phi_x(y_i)\leq0$.
The maximum term is a maximum of affine functions, while $|y-x|$ is affine separately on $[0,x]$ and $[x,\infty)$. Thus $\phi_x$ is convex on each of these two intervals.
Partition the simplex
\[
\bigl\{\mathbf y\in[0,\infty)^{n-1}:\textstyle\sum_i y_i=1-x\bigr\}
\]
into finitely many compact polytopes according to the interval containing each coordinate. The function $\sum_i\phi_x(y_i)$ is convex on every nonempty polytope, so Lemma~\ref{lem:vertex} reduces the proof to its extreme points.

At such an extreme point, at most one coordinate lies outside $\{0,x\}$.
Indeed, if two coordinates $y_i,y_j$ lie in the interiors of their respective intervals, then, for sufficiently small $\varepsilon>0$, both points
\[
\mathbf y^\pm=\mathbf y\pm\varepsilon(\mathbf e_i-\mathbf e_j)
\]
belong to the same polytope, and $\mathbf y=(\mathbf y^++\mathbf y^-)/2$, contrary to extremality.
When $x=0$, the same argument shows that at most one coordinate is positive.
Thus, after permuting the leaves, it suffices to consider
\[
f=(x,\underbrace{x,\ldots,x}_{k\text{ entries}},1-(k+1)x,0,\ldots,0),
\quad 0\leq k\leq n-2,
\]
where
\[
0\leq x\leq u_k:=\min\{c,1/(k+1)\}.
\]
The two variations are
\begin{align*}
\Var_1(M_{\alpha,S_n}f)
&=k(2bx-c)_++\max\{1-(k+1)x-c,b(1-kx)-c,0\},\\
\Var_1(f)&=(n-k-2)x+|1-(k+2)x|.
\end{align*}
The first expression is convex in $x$, and the second is affine on each side of $x=1/(k+2)$.
Their difference in \eqref{eq:sgoal} is therefore convex on each corresponding subinterval. It suffices to check $0$, $u_k$, and the dividing point $1/(k+2)$ whenever it belongs to $[0,u_k]$.
These points have one of three forms: $x=0$, $x=c$, or $x=1/t\leq c$ for an integer $2\leq t\leq n$.
If $u_k=1/(k+1)$, then $c<1$ forces $k\geq1$.

At $x=0$, the function is a single-leaf indicator, and equality holds in \eqref{eq:sgoal}.

At $x=c$, set $\zeta=1-(k+2)c$. The remaining leaf value is $c+\zeta\geq0$, and
\begin{align*}
\Var_1(f)&=nc-1+\zeta+|\zeta|,\\
\Var_1(M_{\alpha,S_n}f)
&=kc(2b-1)+\max\{\zeta,b\zeta+(2b-1)c,0\}.
\end{align*}
If $\zeta\geq0$, then $c\leq1/(k+2)\leq1/2$.
Using $kc=1-2c-\zeta$, the second expression is the maximum of
\[
N_1=(2b-1)(1-2c)+2(1-b)\zeta,
\qquad N_2=(2b-1)(1-c)+(1-b)\zeta.
\]
Subtracting these expressions separately gives
\begin{align*}
(1-c)\Var_1(f)-N_1
&=(1-c)(nc-2b)+(2b-1)c+2(b-c)\zeta\geq0,\\
(1-c)\Var_1(f)-N_2
&=(1-c)(nc-2b)+(1-2c+b)\zeta\geq0.
\end{align*}
Here we used $nc>2b$, $b>c$, and $c\leq1/2$.
If $\zeta\leq0$, then $\Var_1(f)=nc-1$, and
\[
1/c-2\leq k\leq1/c-1,\qquad
\Var_1(M_{\alpha,S_n}f)=\max\{kc(2b-1),b-c-kc(1-b)\}.
\]
The right-hand side is convex as a function of the real variable $k$. At both endpoints, its value is $(2b-1)(1-c)$.
Convexity and $2b-1<nc-1$ now give
\[
\Var_1(M_{\alpha,S_n}f)\leq(2b-1)(1-c)
\leq(1-c)(nc-1)=(1-c)\Var_1(f).
\]

At $x=1/t\leq c$, the center and $t-1$ leaves have value $1/t$, and the remaining $n-t$ leaves have value zero.
Hence
\[
\Var_1(M_{\alpha,S_n}f)=\frac{t-1}{t}(2b-tc)_+,
\qquad \Var_1(f)=\frac{n-t}{t}.
\]
If $tc\geq2b$, the conclusion follows immediately.
Suppose that $tc<2b$, and set $\delta=tc-1$. Since $1/t\leq c$, we have $0\leq\delta<2b-1$.
Substituting $t=(1+\delta)/c$ and rearranging gives the identity
\begin{equation}\label{eq:sflat}
(1-c)(n-t)-(t-1)(2b-tc)
=\frac{(1-c)(nc-2b)-\delta(2b-1-\delta)}{c}.
\end{equation}
If $c\leq1/2$, then \eqref{eq:srel} and
$\delta(2b-1-\delta)\leq(2b-1)^2/4$ imply
\[
(1-c)(nc-2b)-\delta(2b-1-\delta)
\geq b(2b-1)-\frac{(2b-1)^2}{4}
=\frac{(2b-1)(2b+1)}4>0.
\]
If $c>1/2$, the inequality $tc<2b<2$ forces $t=2$ or $3$.
For $t=2$, the left-hand side of \eqref{eq:sflat} is
\[
(n-2)-(n-4)c-2b\geq(n-2)-(n-4)-2=0.
\]
For $t=3$, this difference is $n-3+(9-n)c-4b$.
If $n\geq9$, it is at least $n-3+(9-n)-4=2$; if $5\leq n\leq8$, it is at least
\[
n-3+(9-n)/2-4=(n-5)/2\geq0.
\]
Finally, if $n=4$, then $c=b^2>1/2$, so $b>1/\sqrt2>2/3$ and $3c=3b^2>2b$, contradicting the present assumption.
Every required endpoint therefore satisfies \eqref{eq:sgoal}.
Applying convexity first on each one-dimensional subinterval associated with $k$, and then using the extreme-point bound on each polytope for fixed $x$, proves the assertion for $0\leq x\leq c$.

\paragraph{Case 2: $c\leq x\leq1$.}
The maximal function now equals $x$ at the center. Set
\[
\pi_f=\sum_i(y_i-x)_+,\qquad z_i=\min\{y_i,x\}.
\]
The identities $y_i=z_i+(y_i-x)_+$ and $|y_i-x|=2(y_i-x)_+-(y_i-x)$ give
\[
\sum_i z_i=1-x-\pi_f,\qquad \Var_1(f)=nx-1+2\pi_f.
\]
If $\pi_f>0$, some $y_i>x$, so $1\geq x+y_i>2x$. In particular, $c\leq x<1/2$.
Keeping the original center value $x$ and whole-graph term $c$ fixed, define
\[
L(y)=\max\{y,b(x+y),c\}.
\]
Each affine term has slope in $[0,1]$, so
\[
0\leq L(y')-L(y)\leq y'-y\qquad(y'\geq y).
\]
Also, $L(x)=2bx>x$. For $0\leq z\leq x$,
we have $b(x+z)-z\geq(2b-1)x\geq0$. Thus, setting
\[
e(z)=|\max\{b(x+z),c\}-x|,
\]
yields the termwise estimate
\[
|L(y_i)-x|\leq e(z_i)+(y_i-x)_+.
\]
Equality holds when $y_i\leq x$. When $y_i>x$, the estimate follows from the preceding increment bound and $L(y_i)\geq L(x)>x$.
Consequently,
\[
\Var_1(M_{\alpha,S_n}f)\leq\sum_i e(z_i)+\pi_f.
\]

First suppose that $c\leq x\leq c/b$. Write
\[
e_0=x-c,\qquad \xi=(2b-1)x-e_0=c-2(1-b)x.
\]
Since $\max\{b(x+z),c\}$ is nondecreasing and takes values between $c$ and $2bx$,
$e(z)\leq\max\{e_0,(2b-1)x\}$.
If $\xi\leq0$, then $e(z)\leq e_0$, and hence
\begin{align*}
(1-c)\Var_1(f)-\Var_1(M_{\alpha,S_n}f)
&\geq(1-c)(nx-1+2\pi_f)-(n-1)(x-c)-\pi_f\\
&=(nc-1)(1-x)+(1-2c)\pi_f\geq0.
\end{align*}
The last step uses $nc>1$ and the alternative $\pi_f=0$ or $c<1/2$.
If $\xi\geq0$, then
\[
e(z)\leq e_0+\xi z/x\qquad(0\leq z\leq x).
\]
Indeed, when $\max\{b(x+z),c\}\leq x$, the left-hand side is at most $e_0$; otherwise it equals $b(x+z)-x$, and
\[
e_0+\xi z/x-[b(x+z)-x]=((2-b)x-c)(1-z/x)\geq0.
\]
Summing over the leaves gives
\begin{align*}
&(1-c)\Var_1(f)-\Var_1(M_{\alpha,S_n}f)\\
&\quad\geq(1-c)(nx-1+2\pi_f)-(n-1)e_0
       -\frac{\xi}{x}(1-x-\pi_f)-\pi_f\\
&\quad=(1-x)(nc-2b+1-c/x)+(1-2c+\xi/x)\pi_f\geq0.
\end{align*}
The first term is nonnegative because $nc>2b$ and $x\geq c$. The second vanishes when $\pi_f=0$ and is nonnegative when $\pi_f>0$, since then $c<1/2$ and $\xi\geq0$.

Now suppose that $c/b\leq x\leq1$. Then $bx\geq c$, so $e(z)=|bz-(1-b)x|$ is convex.
The chord joining its endpoint values gives
\[
e(z)\leq(1-z/x)(1-b)x+(z/x)(2b-1)x
=(1-b)x+(3b-2)z.
\]
Using $\sum_i z_i=1-x-\pi_f$, we obtain
\begin{align*}
(1-c)\Var_1(f)-\Var_1(M_{\alpha,S_n}f)
&\geq F(x)+(3b-1-2c)\pi_f,\\
F(x)&=((n+2)b-nc-1)x+1+c-3b.
\end{align*}
The function $F$ is affine, with endpoint values
\[
F(c/b)=(1-c/b)(nc+1-3b),\qquad F(1)=(n-1)(b-c).
\]
By \eqref{eq:srel},
\[
nc+1-3b\geq4b^2-3b+1=4(b-3/8)^2+7/16>0.
\]
Since $c<b$, both endpoint values are nonnegative. Hence $F(x)\geq0$ throughout the interval.
Moreover,
\[
3b-1-2c\geq3b-1-2b^2=(2b-1)(1-b)>0.
\]
Thus \eqref{eq:sgoal} holds in the second case as well.

The two cases cover all normalized nonnegative functions. For a general nonzero nonnegative function, apply the result to $f/\sigma_f$ and restore the original scale by homogeneity.
Lemma~\ref{lem:reduce} extends the bound to all real-valued functions.
The single-leaf indicator considered at the start of the section attains $1-c$, so the constant is sharp.
\end{proof}

\subsection{The exact parameter range for \texorpdfstring{$p=2$}{p=2}}
For the sharp $p=2$ constant of the ordinary maximal operator on star graphs, see \cite{GRM2}*{Theorem 15}.
In the fractional case, we need a further comparison between the coefficients of two-point balls and whole-graph balls. We first prove an auxiliary inequality.
\begin{lem}[A $p=2$ comparison for the two-point-ball coefficient]\label{lem:star-pair-factor}
Let $n\geq4$, $r=\log_2n$, and $1/2<\beta<1$, and set $\gamma=\beta^r$. Then
\[
(\beta-\gamma)\sqrt{1+\frac4{n-2}}\leq1-\gamma.
\]
\end{lem}
\begin{proof}
Since $r\geq2$, we have $0<\gamma<\beta<1$.
The tangent-line bound for the convex function $t\mapsto t^r$ at $t=1$ gives
\[
\gamma=\beta^r\geq1+r(\beta-1),\qquad
1-\gamma\leq r(1-\beta).
\]
Dividing the latter inequality by the positive number $r(1-\gamma)$ and using
$\beta-\gamma=(1-\gamma)-(1-\beta)$, we obtain
\[
\frac{\beta-\gamma}{1-\gamma}
=1-\frac{1-\beta}{1-\gamma}\leq1-\frac1r.
\]
Therefore, since $n=2^r$, it suffices to prove
\[
\biggl(1-\frac1r\biggr)\sqrt{1+\frac4{2^r-2}}\leq1
\qquad(r\geq2).
\]
Both sides are nonnegative. Squaring and multiplying by $r^2(2^r-2)>0$ shows that this is equivalent to
\[
(r-1)^2(2^r+2)\leq r^2(2^r-2),
\]
that is, $g(r):=2^r(2r-1)-4r^2+4r-2\geq0$.
Direct differentiation gives
\[
g(2)=2,\qquad g'(2)=12\log2-4>0,
\]
and
\[
g''(r)=2^r\bigl[(2r-1)(\log2)^2+4\log2\bigr]-8.
\]
Since $\log2>1/2$, for $r\geq2$ we have
\[
g''(r)>4(3/4+2)-8=3>0.
\]
Thus $g'$ is positive on this interval, and $g(r)\geq g(2)>0$.
Combining these estimates yields
\[
\frac{\beta-\gamma}{1-\gamma}\sqrt{1+\frac4{n-2}}
\leq\biggl(1-\frac1r\biggr)\sqrt{1+\frac4{2^r-2}}\leq1.
\]
Multiplication by $1-\gamma>0$ proves the lemma.
\end{proof}
For the parameters of this section, $b^{\log_2n}=(2^{\alpha-1})^{\log_2n}=c$, so the lemma applies with $(\beta,\gamma)=(b,c)$.

\begin{thm}[A necessary and sufficient condition for the $p=2$ sharp constant to equal $1-c$]\label{thm:fracstar-p2-sharp}
Let $n\geq4$, $0<\alpha<1$, and $c=n^{\alpha-1}$. Then
\begin{equation}\label{eq:fs5}
\C_{\alpha,S_n,2}=1-c
\quad\Longleftrightarrow\quad c\geq\frac12
\quad\Longleftrightarrow\quad \alpha\geq1-\log_n2.
\end{equation}
In this range, every single-leaf indicator attains equality. If $c<1/2$, then $\C_{\alpha,S_n,2}>1-c$.
\end{thm}
\begin{proof}
First assume $c\geq1/2$. By Lemma~\ref{lem:reduce} and homogeneity, it suffices for the upper bound to consider $f\geq0$ with $\sigma_f=1$.
Define
\[
Q(f)=\sum_i|M_{\alpha,S_n}f(i)-M_{\alpha,S_n}f(0)|^2
      -(1-c)^2\sum_i(y_i-x)^2.
\]
Since variations are nonnegative, the desired $p=2$ inequality is equivalent to $Q(f)\leq0$.

\paragraph{Sufficiency, Case 1: $0\leq x\leq c$.}
Here $M_{\alpha,S_n}f(0)=c$. Suppose that $Q(f)>0$ for some function.
By continuity, $Q$ attains a positive maximum on the compact set
\[
\bigl\{(x,\mathbf y):0\leq x\leq c,\ y_i\geq0,\ x+\textstyle\sum_i y_i=1\bigr\}
\]
Fix the center value $x$ of a maximizer, and define
\[
\psi_x(y)=\bigl(\max\{y,b(x+y),c\}-c\bigr)^2-(1-c)^2(y-x)^2.
\]
Then $Q(f)=\sum_i\psi_x(y_i)$.
The maximal value at a leaf is strictly greater than $c$ if and only if
\[
y\in I_x:=(u_x,\infty),\qquad
u_x=\min\{c,c/b-x\}>0.
\]
Here $c/b-x\geq c(1/b-1)>0$.

We first show that $\psi_x$ is strictly convex on $I_x$.
The singleton and two-point-ball terms agree at $y_*=bx/(1-b)$. On the corresponding open subintervals of $I_x$,
\[
\psi_x''(y)=
\begin{cases}
2\bigl(b^2-(1-c)^2\bigr),&y<y_*,\\
2\bigl(1-(1-c)^2\bigr),&y>y_*.
\end{cases}
\]
Both expressions are positive because $1-c\leq1/2<b$.
If $y_*\in I_x$, then $y_*>c$, and
\[
\psi'_{x,+}(y_*)-\psi'_{x,-}(y_*)=2(1-b)(y_*-c)>0.
\]
These positive second derivatives and the upward derivative jump prove strict convexity throughout $I_x$.

If a positive maximizer has two leaves $y_i,y_j\in I_x$, choose
$0<\varepsilon<\min\{y_i-u_x,y_j-u_x\}$, and set
$\mathbf y^\pm=\mathbf y\pm\varepsilon(\mathbf e_i-\mathbf e_j)$.
These perturbations preserve nonnegativity and the total sum, and both altered coordinates remain in $I_x$.
Strict convexity in the two altered coordinates gives
\[
Q(x,\mathbf y)
<\tfrac12\bigl[Q(x,\mathbf y^+)+Q(x,\mathbf y^-)\bigr]
\leq Q(x,\mathbf y),
\]
a contradiction. Thus at a positive maximizer, at most one leaf has maximal value greater than $c$.
If there were no such leaf, then $Q(f)=-(1-c)^2\sum_i(y_i-x)^2\leq0$, again a contradiction.

After permuting the leaves, we may therefore assume $M_{\alpha,S_n}f(1)>c$ and denote this leaf value by $y$. All other leaf values belong to $[0,u_x]$.
Set $z=(n-2)^{-1}\sum_{i=2}^{n-1}y_i$, and replace those leaf values by $z$.
Since $z\in[0,u_x]$, the edge differences of the maximal function remain unchanged.
At the same time,
\[
\sum_{i=2}^{n-1}(y_i-x)^2
=(n-2)(z-x)^2+\sum_{i=2}^{n-1}(y_i-z)^2
\geq(n-2)(z-x)^2.
\]
Thus $Q(x,y,z,\ldots,z)\geq Q(f)>0$.
It now suffices to prove the opposite inequality for functions of this form.
Write
\[
\xi=y-x,\qquad \eta=z-x,\qquad
\Var_2(f)^2=\xi^2+(n-2)\eta^2.
\]

If the maximal value at leaf $1$ is attained by the singleton ball, its unique positive edge difference is
\[
y-c\bigl(x+y+(n-2)z\bigr)
=(1-c)\xi+(1-2c)x-c(n-2)z\leq(1-c)\xi.
\]
Positivity forces $\xi>0$, so its square is at most $(1-c)^2[\xi^2+(n-2)\eta^2]$.
If the maximal value is attained by the two-point ball, denote the positive edge difference by $\varepsilon_b$. Then
\begin{align*}
\varepsilon_b
&=b(x+y)-c\bigl(x+y+(n-2)z\bigr)\\
&=(2b-nc)x+(b-c)\xi-c(n-2)\eta.
\end{align*}
When $\eta\geq0$, the inequality $nc>2b$ implies $\varepsilon_b\leq(b-c)\xi$.
Hence $\xi>0$, and $b-c<1-c$ gives the required square bound again.
When $\eta<0$, the condition $z=x+\eta\geq0$ gives $x\geq-\eta$.
Since $2b-nc<0$,
\[
0<\varepsilon_b\leq-(2b-nc)\eta+(b-c)\xi-c(n-2)\eta
=(b-c)(\xi-2\eta).
\]
The Cauchy--Schwarz inequality and Lemma~\ref{lem:star-pair-factor} yield
\begin{align*}
\varepsilon_b
&\leq(b-c)|\xi-2\eta|\\
&\leq(b-c)\sqrt{1+\frac4{n-2}}\sqrt{\xi^2+(n-2)\eta^2}\\
&\leq(1-c)\sqrt{\xi^2+(n-2)\eta^2}.
\end{align*}
Squaring gives $Q(f)\leq0$, contradicting the existence of a positive maximum.
Thus the upper bound holds for all functions with $x\leq c$.

\paragraph{Sufficiency, Case 2: $c\leq x\leq1$.}
Since $c\geq1/2$ and the total sum is one, $y_i\leq1-x\leq x$.
Therefore
\[
b(x+y_i)-y_i\geq(2b-1)x\geq0.
\]
The singleton term at each leaf may thus be omitted, and $M_{\alpha,S_n}f(i)=\max\{b(x+y_i),c\}$.
If $b(x+y_i)<c\leq x$, replacing it by $c$ decreases its distance from $x$; otherwise the maximum is unchanged. Hence
\[
|M_{\alpha,S_n}f(i)-x|^2\leq[b(x+y_i)-x]^2.
\]
Expanding the difference term by term gives
\begin{align*}
Q(f)\leq{}&\bigl(b^2-(1-c)^2\bigr)\sum_i y_i^2
 +2\bigl[(1-c)^2-b(1-b)\bigr]x\sum_i y_i\\
&+(n-1)\bigl[(1-b)^2-(1-c)^2\bigr]x^2.
\end{align*}
The coefficient of the first term is positive. Since $\sum_i y_i=1-x$ and $\sum_i y_i^2\leq(1-x)^2$, the preceding expression is at most
\[
F(x)=(b-x)^2+(n-2)(1-b)^2x^2
-(1-c)^2\bigl[(2x-1)^2+(n-2)x^2\bigr].
\]
Writing $F(x)=A_0+A_1x+A_2x^2$, we have
\begin{align*}
A_0&=b^2-(1-c)^2,\\
A_1&=4(1-c)^2-2b<0,\\
A_2&=1+(n-2)(1-b)^2-(n+2)(1-c)^2.
\end{align*}
Here $4(1-c)^2\leq1<2b$.
If $A_2\leq0$, then $F'(x)=A_1+2A_2x<0$ for $x\geq0$, so $F(x)\leq F(c)$.
If $A_2>0$, then $F$ is convex, and hence
\[
F(x)\leq\frac{1-x}{1-c}F(c)+\frac{x-c}{1-c}F(1)
\leq\max\{F(c),F(1)\}\quad(c\leq x\leq1).
\]
It therefore suffices to show that both endpoint values are negative.

At $x=1$,
\[
F(1)=(n-1)\bigl[(1-b)^2-(1-c)^2\bigr]<0.
\]
At $x=c$, set $\delta=b-c>0$. Using
$(1-b)^2-(1-c)^2=-\delta(2-b-c)$, we obtain
\[
F(c)=\delta\bigl[\delta-(n-2)c^2(2-b-c)\bigr]
-(1-c)^2(2c-1)^2.
\]
Lemma~\ref{lem:star-pair-factor} gives
\[
\frac{\delta}{1-c}\leq\sqrt{\frac{n-2}{n+2}}
<\frac{2(n-2)}{n+2}.
\]
The strict inequality holds because $n\geq4$.
Rearranging yields
\[
4\delta<(n-2)\bigl[2(1-c)-\delta\bigr]=(n-2)(2-b-c).
\]
Since $c^2\geq1/4$ and $2-b-c>0$, this implies
$\delta<(n-2)c^2(2-b-c)$, and hence $F(c)<0$.
Thus $Q(f)\leq F(x)<0$, completing the second case.

Restoring the original scale and using the reduction to nonnegative functions proves the upper bound for all real-valued functions.
A single-leaf indicator attains $1-c$, so this is the sharp constant when $c\geq1/2$.

\paragraph{Necessity: $c<1/2$.}
Fix $n,\alpha$ satisfying this condition, and take
\[
f_t(0)=t,\qquad f_t(1)=1+t,\qquad
f_t(i)=0\quad(2\leq i\leq n-1),
\]
where
\[
0<t<\min\biggl\{\frac{c}{1-2c},\frac{1-b}{2b-1},
\frac{2(1-2c)}{(n-2)(1-c)}\biggr\}.
\]
The first two restrictions give, respectively,
\[
t<c(1+2t),\qquad b(1+2t)<1+t.
\]
Also, $1+t-c(1+2t)=1-c+(1-2c)t>0$, and $bt<t<c(1+2t)$.
Thus \eqref{eq:f1} gives
\[
M_{\alpha,S_n}f_t(1)=1+t,\qquad
M_{\alpha,S_n}f_t(i)=c(1+2t)\quad(i\ne1).
\]
The original absolute edge differences are $1,t,\ldots,t$, while the maximal function has exactly one nonzero edge difference. Hence
\[
\frac{\Var_2(M_{\alpha,S_n}f_t)}{\Var_2(f_t)}
=\frac{1-c+(1-2c)t}{\sqrt{1+(n-2)t^2}}.
\]
Subtracting $(1-c)^2[1+(n-2)t^2]$ from the square of the numerator gives
\begin{align*}
&2(1-c)(1-2c)t+\bigl[(1-2c)^2-(n-2)(1-c)^2\bigr]t^2\\
&\qquad>2(1-c)(1-2c)t-(n-2)(1-c)^2t^2>0.
\end{align*}
The last step uses the third restriction on $t$.
This variation ratio is therefore strictly greater than $1-c$, so $\C_{\alpha,S_n,2}>1-c$.
Necessity follows. Finally, $n^{\alpha-1}\geq1/2$ is equivalent to $\alpha\geq1-\log_n2$, proving the equivalence of all three conditions.
\end{proof}

\subsection{Counterexamples and parameter restrictions}
Theorem~\ref{thm:fracstar-p2-sharp} determines the exact parameter range of the formula $1-c$ for the fixed exponent $p=2$.
Varying the variation exponent, we now show that the $p=1$ formula does not extend to all $p<1$, and that the $p=2$ formula does not extend to all $p>2$ even under the condition $c\geq1/2$.
\begin{prop}[Two families of test functions for the parameter restrictions]\label{prop:fracstar-counter}
The functions in the following two assertions are nonnegative and have positive variation. The parameter sequences are chosen separately in the two assertions.
\begin{enumerate}[label=\textup{(\roman*)}]
\item For each fixed $0<p<1$, there exist $\alpha_n\downarrow0$ and functions $f_n$ on $S_n$ such that
\[
\frac{\Var_p(M_{\alpha_n,S_n}f_n)}{\Var_p(f_n)}\longrightarrow\infty.
\]
\item There exist arbitrarily large $n$, parameters $\alpha_n\in(0,1)$ satisfying $n^{\alpha_n-1}=1/2$, exponents $p_n>2$, and functions $f_n$ on $S_n$ such that
\[
\frac{\Var_{p_n}(M_{\alpha_n,S_n}f_n)}{\Var_{p_n}(f_n)}>\frac12.
\]
\end{enumerate}
\end{prop}
\begin{proof}
For the first assertion, fix $0<p<1$. For every $n\geq4$, take
\[
\alpha_n=\log_n\frac{n}{n-1},\qquad
f_n=(1,\underbrace{1,\ldots,1}_{n-2\text{ leaves}},0).
\]
Then $0<\alpha_n<1$, and $c=n^{\alpha_n-1}=1/(n-1)$.
The numerator $\log(n/(n-1))$ decreases strictly to zero, while $\log n$ increases strictly. Thus $\alpha_n\downarrow0$.
The original function has only one nonzero edge difference, namely the edge joining the zero leaf to the center, and its absolute value is $1$. Hence $\Var_p(f_n)=1$.
Also, $\sigma_{f_n}=n-1$, and the whole-graph term is $c\sigma_{f_n}=1$. Formula \eqref{eq:f1} gives
\[
M_{\alpha_n,S_n}f_n(0)=M_{\alpha_n,S_n}f_n(n-1)=1,
\qquad
M_{\alpha_n,S_n}f_n(i)=2^{\alpha_n}\quad(1\leq i\leq n-2).
\]
At each leaf of value $1$, the two-point-ball value $2^{\alpha_n}>1$ is maximal. At the zero leaf, the two-point-ball value is $2^{\alpha_n-1}<1$, so the whole-graph term is maximal.
Therefore
\[
\frac{\Var_p(M_{\alpha_n,S_n}f_n)}{\Var_p(f_n)}
=(n-2)^{1/p}(2^{\alpha_n}-1).
\]
Since $\log(1+t)\sim t$ and $e^t-1\sim t$ as $t\to0$,
\[
\alpha_n\sim\frac1{n\log n},\qquad
2^{\alpha_n}-1\sim\frac{\log2}{n\log n}.
\]
The preceding ratio is asymptotic to $(\log2)n^{1/p-1}/\log n$ and tends to infinity because $1/p-1>0$.

For the second assertion, take, for each $n\geq16$,
\[
\alpha_n=1-\log_n2,\qquad p_n=\log_2(n-1),\qquad
f_n=(1,2,0,\ldots,0).
\]
Then $0<\alpha_n<1$, $p_n>2$, and
\[
c=n^{\alpha_n-1}=\frac12,\qquad
(n-1)^{1/p_n}=2,\qquad
b=2^{\alpha_n-1}\geq2^{-1/4}>\frac56.
\]
Since $\sigma_{f_n}=3$, the maximal value at the center is $3/2$. At a zero leaf, the three ball values are $0,b,3/2$, so its maximal value is also $3/2$. At leaf $1$, they are $2,3b,3/2$, with maximum $3b$.
Consequently,
\[
\Var_{p_n}(f_n)=(n-1)^{1/p_n}=2,\qquad
\Var_{p_n}(M_{\alpha_n,S_n}f_n)=3b-\frac32,
\]
and hence
\[
\frac{\Var_{p_n}(M_{\alpha_n,S_n}f_n)}{\Var_{p_n}(f_n)}
=\frac32\biggl(b-\frac12\biggr)>\frac12.
\]
This proves both assertions.
\end{proof}

\subsection{The sharp infinite-exponent constant}
For the infinite exponent, it suffices to control the absolute difference along each edge.
Monotone changes in the other leaf values reduce the comparison of any positive directed edge difference to an explicit one-parameter family.
\begin{thm}[The sharp infinite-exponent constant on star graphs]\label{thm:fracstarinf}
Let $n\geq4$, $0<\alpha<1$, $b=2^{\alpha-1}$, and $c=n^{\alpha-1}$. Set
\[
A_*=\min\{1,1-c+(1-2c)_+\},\qquad
B_*=
\begin{cases}
3(b-c),&c\geq1/3,\\[2pt]
b+\dfrac{(2b-1)c(n-3)}{nc-1},&c<1/3.
\end{cases}
\]
Then
\[
\C_{\alpha,S_n,\infty}=\max\{A_*,B_*\}.
\]
An extremizer can be chosen from the family
\[
f_x=(x,x+1,(x-1)_+,\ldots,(x-1)_+),\qquad x\geq0.
\]
\end{thm}
\begin{proof}
By the reduction to nonnegative functions and homogeneity, it suffices to consider $f\geq0$ with $\Var_\infty(f)=1$.
Thus every leaf value satisfies
\[
(x-1)_+\leq y_i\leq x+1.
\]
We estimate the two possible directions of the maximal-function edge difference separately.

\paragraph{The maximal value at the leaf is smaller than at the center.}
If $M_{\alpha,S_n}f(0)>M_{\alpha,S_n}f(i)$, then, since the leaf maximal value is at least the whole-graph term, the maximum at the center must be attained by the singleton ball. Thus $M_{\alpha,S_n}f(0)=x$.
When $0\leq x\leq1$,
\[
M_{\alpha,S_n}f(0)-M_{\alpha,S_n}f(i)
\leq x-b(x+y_i)\leq(1-b)x\leq1-b.
\]
When $x\geq1$, the inequalities $y_i\geq x-1$ and $1-2b<0$ give
\[
M_{\alpha,S_n}f(0)-M_{\alpha,S_n}f(i)
\leq x-b(2x-1)=(1-2b)x+b\leq1-b.
\]
Hence every difference in this direction is at most $1-b$.

\paragraph{The maximal value at the leaf is larger than at the center.}
Fix the leaf value $y$ under consideration and the center value $x$.
Since the whole-graph term does not exceed the maximal value at the center, the positive part of the difference in this direction is
\[
\max\{0,y-M_{\alpha,S_n}f(0),b(x+y)-M_{\alpha,S_n}f(0)\}.
\]
Lowering the other $n-2$ leaf values to $z=(x-1)_+$ decreases the whole-graph term and the maximal value at the center, and therefore cannot decrease this expression.
The differences associated with the singleton and two-point balls are then, respectively,
\[
\min\{y-x,(1-c)y-c(x+(n-2)z)\},
\]
\[
\min\{(b-1)x+by,(b-c)x+(b-c)y-c(n-2)z\}.
\]
Here we used $a-\max\{u,v\}=\min\{a-u,a-v\}$.
The slopes in $y$ of the affine terms are $1,1-c,b,b-c$, all positive. Thus both minimum functions are nondecreasing.
Raising $y$ to $x+1$ therefore cannot decrease the positive difference under consideration.
The resulting function is precisely $f_x$. Every absolute edge difference is at most $1$, and the difference along $\{0,1\}$ is exactly $1$. Hence $\Var_\infty(f_x)=1$.

For this family,
\[
\sigma_{f_x}=
\begin{cases}
2x+1,&0\leq x\leq1,\\
nx-(n-3),&x\geq1.
\end{cases}
\]
Define
\[
H_1(x)=x+1-\max\{x,c\sigma_{f_x}\},\qquad
H_2(x)=b(2x+1)-\max\{x,c\sigma_{f_x}\}.
\]
The positive part of the difference at the distinguished leaf is $\max\{0,H_1(x),H_2(x)\}$.
The preceding edgewise comparisons give
\[
\Var_\infty(M_{\alpha,S_n}f)
\leq\max\bigl\{1-b,\sup_{x\geq0}H_1(x),\sup_{x\geq0}H_2(x)\bigr\}.
\]
Conversely, for every $x\geq0$,
\[
\Var_\infty(M_{\alpha,S_n}f_x)\geq\max\{0,H_1(x),H_2(x)\},
\qquad \Var_\infty(f_x)=1.
\]
Since $H_1(0)=1-c>1-b$, it remains to maximize $H_1,H_2$ and verify attainment.

\paragraph{The function associated with the singleton ball.}
Substitution of the formula for the total sum gives
\[
H_1(x)=
\begin{cases}
\min\{1,1-c+(1-2c)x\},&0\leq x\leq1,\\
\min\{1,1+c(n-3)-(nc-1)x\},&x\geq1.
\end{cases}
\]
Since $nc>1$, the second branch is nonincreasing.
The first branch is nondecreasing when $c\leq1/2$, with maximum $H_1(1)=\min\{1,2-3c\}$. When $c\geq1/2$, it is nonincreasing, with maximum $H_1(0)=1-c$.
The branches agree at $x=1$, so
\[
\max_{x\geq0}H_1(x)=\min\{1,1-c+(1-2c)_+\}=A_*.
\]
In particular, $A_*\geq1-c>0$.

\paragraph{The function associated with the two-point ball.}
For $0\leq x\leq1$,
\[
H_2(x)=\min\{(2b-1)x+b,2(b-c)x+b-c\},
\]
and both terms are strictly increasing.
For $x\geq1$, the whole-graph term at the center minus its singleton value is
\[
c\sigma_{f_x}-x=(nc-1)x-c(n-3).
\]
If $c\geq1/3$, this difference is already $3c-1\geq0$ at $x=1$, and its slope $nc-1$ is positive.
Thus the whole-graph ball always attains the maximum at the center when $x\geq1$, and the slope of $H_2$ on this interval is $2b-nc<0$.
Together with monotonicity on the preceding interval, this shows that the maximum is attained at $x=1$, with value
\[
\max_{x\geq0}H_2(x)=H_2(1)=3(b-c).
\]
If $c<1/3$, the two ball values at the center agree at
\[
x_*:=\frac{c(n-3)}{nc-1}>1
\]
where the last inequality is equivalent to $1-3c>0$.
For $1\leq x\leq x_*$, the singleton ball attains the maximum at the center, so $H_2(x)=(2b-1)x+b$ has positive slope.
For $x\geq x_*$, the whole-graph ball attains the maximum, and the slope of $H_2$ is $2b-nc<0$.
The two expressions agree at $x_*$, so the maximum is
\[
H_2(x_*)=b+\frac{(2b-1)c(n-3)}{nc-1}.
\]
This proves $\max_{x\geq0}H_2(x)=B_*$.
At $c=1/3$, the latter expression also gives $b+(2b-1)=3(b-c)$, so the two formulas agree.

The preceding bounds show that every normalized function satisfies
\[
\Var_\infty(M_{\alpha,S_n}f)\leq\max\{A_*,B_*\}.
\]
To verify attainment, suppose first that $A_*\geq B_*$. Choose $f_0$ if $c\geq1/2$ and $f_1$ if $c<1/2$. If $B_*>A_*$, choose $f_1$ when $c\geq1/3$ and $f_{x_*}$ when $c<1/3$.
Each chosen function has infinite-exponent variation $1$, and the positive difference of its maximal function along $\{0,1\}$ equals $\max\{A_*,B_*\}$.
Thus the upper bound is attained, proving the sharp-constant formula.
\end{proof}

\subsection{Sharp mixed constants}
For a fixed finite vector of edge differences, the unnormalized $\ell^r$ norm or quasi-norm is nonincreasing in the exponent.
Applying this fact separately to the original function and its maximal function yields the following mixed estimates from the bounds for $p=1$ and $p=2$.
\begin{thm}[Sharp mixed constants on star graphs]\label{thm:fracstar-mixed-sharp}
Let $n\geq4$, $0<\alpha<1$, and $c=n^{\alpha-1}$.
\begin{enumerate}[label=\textup{(\roman*)}]
\item If $0<p\leq1\leq q\leq\infty$, then $\C_{\alpha,S_n;p,q}=1-c$.
\item If $c\geq1/2$ and $0<p\leq2\leq q\leq\infty$, then $\C_{\alpha,S_n;p,q}=1-c$.
\end{enumerate}
In both cases, every single-leaf indicator attains equality.
\end{thm}
\begin{proof}
If $0<p\leq1\leq q\leq\infty$, Theorem~\ref{thm:fracstar-p1-sharp} and monotonicity of finite-sequence norms give
\[
\Var_q(M_{\alpha,S_n}f)\leq\Var_1(M_{\alpha,S_n}f)
\leq(1-c)\Var_1(f)\leq(1-c)\Var_p(f).
\]
If $c\geq1/2$ and $0<p\leq2\leq q\leq\infty$, Theorem~\ref{thm:fracstar-p2-sharp} gives, in the same way,
\[
\Var_q(M_{\alpha,S_n}f)\leq\Var_2(M_{\alpha,S_n}f)
\leq(1-c)\Var_2(f)\leq(1-c)\Var_p(f).
\]
A single-leaf indicator has one nonzero edge difference before and after application of the maximal operator, of absolute values $1$ and $1-c$, respectively.
It attains equality for every stated exponent pair.
\end{proof}

\begingroup
\raggedright

\endgroup

\begin{bibdiv}
\begin{biblist}
\bib{GKM}{article}{
  author={Gonz\'alez-Riquelme, Cristian},
  author={Kova\v{c}, Vjekoslav},
  author={Madrid, Jos\'e},
  title={Sharp variational inequalities for the Hardy--Littlewood maximal operator on finite undirected graphs},
  date={2026},
  eprint={2603.12462v1},
}
\bib{GRM1}{article}{
  author={Gonz\'alez-Riquelme, Cristian},
  author={Madrid, Jos\'e},
  title={Sharp inequalities for maximal operators on finite graphs},
  journal={J. Geom. Anal.},
  volume={31},
  number={10},
  date={2021},
  pages={9708--9744},
  doi={10.1007/s12220-021-00625-0},
}
\bib{GRM2}{article}{
  author={Gonz\'alez-Riquelme, Cristian},
  author={Madrid, Jos\'e},
  title={Sharp inequalities for maximal operators on finite graphs, II},
  journal={J. Math. Anal. Appl.},
  volume={506},
  number={2},
  date={2022},
  pages={Paper No. 125647},
  doi={10.1016/j.jmaa.2021.125647},
}
\bib{Kinnunen}{article}{
  author={Kinnunen, Juha},
  title={The Hardy--Littlewood maximal function of a Sobolev function},
  journal={Israel J. Math.},
  volume={100},
  date={1997},
  pages={117--124},
  doi={10.1007/BF02773636},
}
\bib{LLW}{article}{
  author={Liu, Feng},
  author={Lian, Yongzhi},
  author={Wang, Yanzhen},
  title={Sharp $L^p$ and variation inequalities for maximal operators on finite connected graphs},
  date={2026},
  eprint={2611.15377v1},
}
\bib{LX}{article}{
  author={Liu, Feng},
  author={Xue, Qingying},
  title={On the variation of the Hardy--Littlewood maximal functions on finite graphs},
  journal={Collect. Math.},
  volume={72},
  date={2021},
  pages={333--349},
  doi={10.1007/s13348-020-00290-6},
}
\bib{Ore}{article}{
  author={Ore, O.},
  title={Diameters in graphs},
  journal={J. Combin. Theory},
  volume={5},
  number={1},
  date={1968},
  pages={75--81},
  doi={10.1016/S0021-9800(68)80030-4},
}
\end{biblist}
\end{bibdiv}
\end{document}